\documentclass{article}
\usepackage{amsmath,amsthm,amsfonts,amssymb}

\usepackage[dvips]{graphicx}
\usepackage[usenames,dvipsnames]{color}
\usepackage{fancyhdr}

\newtheorem{theorem}{Theorem}[section]
\newtheorem{remark}[theorem]{Remark}
\newtheorem{lemma}[theorem]{Lemma}

\newcommand{\nn}{\nonumber}

\theoremstyle{definition}

\newcommand {\p}{\partial}
\newcommand{\e}{\varepsilon}
\newcommand{\bor}{\mbox{\boldmath $\varrho$}}
\newcommand{\mbb}[1] {\mathbb{#1}}
\newcommand{\mb}[1] {\mathbf{#1}}
\newcommand{\mc}[1] {\mathcal{#1}}
\def\newline{\ifvmode \else\unskip\hfil\penalty -10000\relax\fi}

\begin{document}

\centerline{\Large{\bf{On a macroscopic limit}}}
\centerline{\Large{\bf{of a kinetic model of alignment}}}

\bigskip\bigskip
\begin{center}
\textbf{Jacek Banasiak}
\smallskip

School of Mathematical  Sciences, University of KwaZulu-Natal \\
Durban, South Africa \\
Institute of Mathematics, Technical University of \L\'{o}d\'{z},\\
\L\'{o}d\'{z}, Poland\\
\texttt{banasiak@ukzn.ac.za} \\

\bigskip
and

\bigskip
\textbf{Miros{\l}aw Lachowicz}
\smallskip

Institute of Applied Mathematics and Mechanics \\
Faculty of Mathematics, Informatics and Mechanics, University of Warsaw \\
ul. Banacha 2, 02--097 Warsaw, Poland \\
\texttt{lachowic@mimuw.edu.pl}
\end{center}

\bigskip
\noindent
{\small{\textbf{Abstract}.
In the present paper the macroscopic limits of the kinetic model for interacting entities
(individuals, organisms, cells) proposed in \cite{ADL} are studied.
The kinetic model is one--dimen\-sio\-nal and entities are characterized
by their position and orientation (+/-) with swarming interaction
controlled by the sensitivity parameter $\gamma$. The macroscopic
limits of the model are considered for solutions close either to
the diffusive (isotropic) or to the aligned (swarming) equilibrium
states for various $\gamma$. In the former case the classical linear
diffusion equation results whereas in the latter a traveling wave
solution does both in the zeroth (`Euler') and first (`Navier--Stokes')
order of approximation.}}

\bigskip\bigskip
\noindent
\section{Introduction}\label{intro}

\medskip
The paper \cite{ADL} focused on biological alignment (swarming) characterized as coherent
motion of groups of entities (individuals, organisms or cells) into the same direction.
We understand alignment as an adaptation of orientation to that of neighborhood.
It results from the interplay of the behaviors of entities by local interactions.
Swarms, also called herds, flocks, schools, clusters depending on whether they refer,
respectively, to insects, mammals, birds, fish,  bacteria and cells are often observed
in nature. The typical examples are herds of sheep, flocks of birds or schools of fish
(see References in \cite{ADL} and \cite{DD, MS, WJKA, Z}). In some of these phenomena
alignment is a typical behavior.

The phenomena of alignment (swarming), viewed as a cooperative phenomenon arising from
the interaction of a number of entities, can be analyzed by mathematical models. The entities 
are characterized by their position, orientation (velocity) and migrate in space and time.
Various microscopic, mesoscopic and macroscopic models describing swarming have been proposed. The dynamics of
alignment in the space homogeneous case was considered e.g. in \cite{GLM, ME}. Microscopic models,
for example the cellular automaton approach, \cite{BDG, DD}, and a simulation model by Ben--Jacob
{\textsl{et al}}. \cite{BS} allows to  distinguish individual entities.  By means of a macroscopic
model based e.g. on hydrodynamics (see \cite{TT}), it is possible to study the dynamics
of densities of entities. The review of various individual--based models based on ODEs is
presented in \cite{CCR}. In a number of papers the authors have studied stochastic models for swarm aggregation, where the individuals, in addition to the classical
Brownian random dispersal, are subject to long range attraction and short range repulsion. For instance, see \cite{BCM} and references therein, under suitable laws of large numbers it was shown that, for a large number
of individuals, the evolution of the empirical distribution of the population can be approximated using a suitable limit nonlinear degenerate and nonlocal parabolic equation. In \cite{EH} a class of stochastic
individual-based models, written in terms of coupled velocity jump processes, was presented and analyzed.
Moreover, in the limit of large populations, a system of two kinetic equations with nonlocal and nonlinear
right hand side was derived and analyzed.

Mesoscopic models have been considered e.g. in \cite{Lu, ADL, PSV, CFRT, AIR}; see
references therein. Such models are usually of Boltzmann type, i.e., they are related to
a statistical description of one test entity.

The relations between the microscopic, mesoscopic and macroscopic models is one of the most
challenging mathematical problem. Its importance is particularly visible in justification of the
macroscopic models by a well-defined behavior of the entities of the system -- see e.g.
\cite{BCM, CL, CCR, DM, EH, La, La2, MBbook, MCO} and references therein.

In the present paper we study the macroscopic limits for the discrete velocity Boltzmann--type
(mesoscopic) model based on a well-defined microscopic `\textsl{majority-choice}' interaction
introduced in \cite{ADL}. We are interested in the formation of alignment of the population of
entities that do not possess a leader but polarization arises as a result of local alignment
interactions.

The model, proposed in \cite{ADL}, is one--dimen\-sio\-nal and the entities are characterized
by their position and orientation (+/-) while the alignment (swarming) interaction is
characterized by a sensitivity parameter $\gamma$.
In the present paper the macroscopic limits of the model are
considered for solutions either close to the diffusive (isotropic), or to the aligned (swarming),
equilibrium states for various values of $\gamma$. In the former case, the classical linear
diffusion equation results whereas in the latter we obtain a  traveling wave solution both in the
zeroth (`Euler') and the first (`Navier-Stokes') order of approximation. An interesting feature
of the considered model is that the nonlinear kinetic model in the `{\textsl{hydrodynamic}}'
limits gives, in both cases, linear macroscopic equations.
Moreover, it follows, at least formally, that all terms of the asymptotic expansion of the
bulk part of the approximation, which are complementary to the traveling wave solution (its'
`{\textsl{kinetic}}' part), vanish and thus the bulk part of the approximation only consists of the
traveling wave, whereas the complementary part is made up only of the initial layer  terms.

The plan of the paper is as follows. Section \ref{model} presents the model that was introduced in
\cite{ADL}. In Section \ref{existence} we prove a unified existence and uniqueness result valid
both in weak and strong interaction cases. Section \ref{formal} deals with formal macroscopic
limits of the model. This section also indicates some interesting open problems: the hydrodynamic
limit in the `{\textsl{diffusive}}' case. Finally we propose an asymptotic result in the
`{\textsl{aligned}}' case. The result is proved in the spirit of
the Tikhonov--Vasileva theory (cf. \cite{TVS}) that
was formulated only in the finite--dimensional dynamical system case. Therefore it is an example
of a generalization of the theory available for ODEs into some infinite dimensional case.

\bigskip
\section{The model}\label{model}

\medskip
Let $f=f(t,j,x)$ be the probability that an individual is at time $t$ at
point $x$ with orientation $j$, where $t>0$, $x\in{\mathbb{R}}$ and
$j\in\{ -1,1\}$. Here the orientation plays the role of discrete velocity in the sense that it give the direction of motion.

We describe the migration together with the changes of the orientation variable where the latter is modeled
by an  interaction operator. We have the following general expression (valid for
every migration-interaction case):
\begin{equation}\label{2.1}
f (t+\Delta t,j,x+j\Delta x) = f (t,-j,x)\,\mathfrak{P} (t,-j,x)  +
f (t,j,x) \,\mathfrak{P}' (t,j,x) \,,
\end{equation}
where
\[
\mathfrak{P}(t,j,x) = {\mathrm{Probability}}\Big( \mathrm{a}\; {\mathrm{change}}\;{\mathrm{of}}\;
{\mathrm{orientation}}\;{\mathrm{in}}\; \Delta t \; \Big|\; (t,j,x) \Big)\,,
\]
is the probability of a change of orientation in the interval of time $[t,t+\Delta t]$
of an individual that is at $t$ at position $x$ with orientation $j$ and
\[
\begin{array}{lll}
& \mathfrak{P}' (t,j,x) = \\
& {\mathrm{Probability}}\Big( {\mathrm{no}}\;{\mathrm{changes}}\;{\mathrm{of}}\;
{\mathrm{orientation}}\;{\mathrm{in}}\; \Delta t \;\Big|\; (t,j,x) \Big) = \\
& 1-{\mathfrak{P}}(t,j,x)\,.
\end{array}
\]

Different choices of the probabilities $\mathfrak{P}$ give rise to different
models, see e.g. \cite{ADL} where, in particular, we  proposed
\begin{equation}\label{2.2}
\mathfrak{P} (t,j,x) = {\frac{ \chi\Big( \sum\limits_{k,l=\pm 1} f (t,k,x+al
)>0 \Big) {\Big( \sum\limits_{l=\pm 1} f
(t,-j,x+al)\Big)^{\gamma}}}{\Big( \sum\limits_{l=\pm 1}
f (t,-j, x+al)\Big)^{\gamma} + \Big(\sum\limits_{l=\pm}
f (t,j, x+al)\Big)^{\gamma}}} \;\;\Delta t \,,
\end{equation}
where $a>0$ and $\gamma >0$ are parameters and
$\sum\limits_{l=\pm 1} f (t,j, x+al)$ is the neighborhood
density in direction $j$, $\chi ({\rm true}) =1$, $\chi ({\rm
false}) =0$. The parameter $\gamma$ describes sensitivity of the
interaction. If $\gamma$ is small (close to $0$), then the
probability of a change of orientation only weakly depends on the
actual orientation. On the other hand, for large $\gamma$  the
probability of a change of orientation strongly depends on the
actual orientation.

Assuming that $a=\Delta x$ and $\Delta x = j\,\Delta
t$, in the limit $\Delta t \to 0$, we obtain the following system
of two equations
\begin{eqnarray}
 &&\partial_t f (t,j,x) +j\partial_x f (t,j,x)
=
{\frac{\chi\big( \sum\limits_{k=\pm 1} f (t,k,x)>0 \big)}
{\big(f (t,-j, x)\big)^{\gamma} + \big( f (t,j,
x)\big)^{\gamma}}}\nn\\
&&\phantom{xxxxx} \times\bigg( f(t,-j,x)\Big(f(t,j,x)\Big)^{\gamma} - f(t,j,x)
\Big( f(t,-j,x)\Big)^{\gamma}\bigg)\,, \quad
 j=-1,1\,.\nn\\
\label{2.3}
\end{eqnarray}

For $\gamma =1$, Eq. (\ref{2.3}) decouples into two free-streaming equations.
Throughout the paper we consider either $\gamma\in ]0,1[$ or
$\gamma >1$.

The non-negative and non-zero  equilibrium solutions $\bar{f},$
corresponding to the space homogeneous version of Eq. (\ref{2.3}), can be found from
\begin{equation}\label{2.4}
\bar{f}_{-1} \Big({\bar{f}}_1  \Big)^{\gamma}
- {\bar{f}}_1 \Big({\bar{f}}_{-1}  \Big)^{\gamma} = 0.
\end{equation}
The only non-zero solutions $\bar f =(\bar f_1,\bar f_{-1})$ of (\ref{2.4})   are given by
\begin{equation}\label{2.5}
{\bar{f}}_1 ={\bar{f}}_{-1}>0\,,
\end{equation}
or
\begin{equation}\label{2.6}
{\bar{f}}_j =0\,\quad\textrm{and}\;\;
{\bar f}_{-j}>0\,,\qquad \textrm{for}\;\textrm{some}\;
j\in\{ -1,1\}.
\end{equation}

Equation (\ref{2.5}) corresponds to equal probabilities of both
orientations ({\textsl{a diffusive picture}}), whereas Eq. (\ref{2.6}) is related
to an {\textsl{aligned picture}}.

In the space-homogeneous case the trajectories of the corresponding
ODEs are contained in the straight lines defined by
\begin{equation}\label{2.7}
f_{-1} + f_1 = c\,,\qquad \mathrm{where}\quad
c=f_{-1} (0)+ f_1 (0) >0\,,
\end{equation}
and have different types of behavior according to the value of $\gamma$.
For $\gamma > 1$
\begin{equation}\label{2.8}
\lim\limits_{t\to\infty}f_{-1} =0\,,\quad
\lim\limits_{t\to\infty}f_1 = f_1 (0)+ f_{-1}
(0) \,,\quad {\mathrm{for}}\; f_1 (0) > f_{-1} (0)\,,
\end{equation}
\begin{equation}\label{2.9}
\lim\limits_{t\to\infty}f_{-1}=f_1 (0) + f_{-1}(0)\,,\quad
\lim\limits_{t\to\infty}f_1 = 0\,,\quad \mathrm{for}\;
f_1 (0) < f_{-1} (0)\,,
\end{equation}
whereas for $0<\gamma <1$ we have
\begin{equation}\label{2.10}
\lim\limits_{t\to\infty}f_1 =
\lim\limits_{t\to\infty}f_{-1} = {\frac{f_1 (0)
+f_{-1} (0) }{2}}\,.
\end{equation}

The important problem whether the similar behavior can be observed in the
spa\-tia\-lly non\-homo\-ge\-neous case remains open.

We note that the spatially nonhomogeneous  Eq.
(\ref{2.3}) admits simple solutions \begin{equation}\label{2.11}
f (t,j,x) =f (t,-j,x) = {\mathrm{const}}.\geq 0 \,,
\end{equation}
and
\begin{equation}\label{2.12}
f (t,j,x) =0\,,\qquad f (t,-j,x) =\phi
(x+j\, t)\,,\quad j=1, -1,
\end{equation}
 where $\phi$ is a given non--negative function, obtained  by equating  both sides of (\ref{2.3})  to zero.

Thus we observe that  the possible asymptotic behavior for Eq. (\ref{2.3}) may differ from that
of the Carleman-type equations studied in \cite{P, GS, ST, SV}.

\section{Existence}\label{existence}
The (global) existence and uniqueness results for $\gamma >1$ as well as the existence
of an {\sl entropy} functional  were established
in \cite{ADL} taking advantage of the Lipschitz continuity and the conservativeness of the system.

Here we propose the global existence and
uniqueness results in an appropriate Banach space for both $\gamma >1$ and
$\gamma < 1$. We will consider one dimensional periodic boundary conditions, that is, we assume
$x\in{\mathbb{T}}$, but similar results are possible for $x\in{\mathbb{T}}^d$, where $d\geq 1$ and ${\mathbb{T}}^d$ is
a $d$-dimensional torus.

Let $X_1$ be the space of real--valued functions equipped with the norm
\[
\|\, f\,\|_1 =\sum\limits_{j=\pm 1} \int\limits_{\mathbb{T}} |f(j,x)|\,{\mathrm{d}}x\,.
\]
The cone of non--negative functions in $X_1$
is denoted by $X_1^+$. We define the following operators
\[
Q[f]= Q^+ [f] - Q^- [f]\,,
\]
where
\[
Q^+ [f](j,x) = \chi\Big(\sum\limits_{k=\pm 1} f(k,x)>0\Big)
{\frac{f^{\gamma}(j,x)f(-j,x)}{f^{\gamma}(j, x) + f^{\gamma}(-j, x)}}\,,
\]
and
\[
Q^- [f](j,x) = \chi\Big(\sum\limits_{k=\pm 1} f(k,x)>0\Big)
{\frac{f(j,x)f^{\gamma} (-j,x)}{f^{\gamma} (j, x) + f^{\gamma} (-j, x)}}\,.
\]
Note that, for $f>0$,
\begin{equation}\label{ba-ba}
Q[f]=R[f]-f\,,
\end{equation}
where
\[
R[f] (j,x) = f^{\gamma} (j,x)
{\frac{f(1,x)+ f(-1,x)}{f^{\gamma} (1,x)+ f^{\gamma} (-1,x)}}\,.
\]

We denote
\[
\mathbb{Q}(h,g)(j,x) = {\frac{h^{\gamma}(j,x)h(-j,x)}{h^{\gamma}(j, x) + g^{\gamma}(-j, x)}}\,,
\]
and
\[
\mathbb{P}(h,g)(j,x) =
{\frac{h^{\gamma} (-j,x)}{g^{\gamma} (j, x) + h^{\gamma} (-j, x)}}\,.
\]

Let
\begin{equation}\label{szarp}
f^{\sharp} (t,j,x) = f(t,j,x+jt)\,.
\end{equation}

We define the following sequences $\{ g_n \}_{n\in\mathbb{N}}$ and $\{ h_n \}_{n\in\mathbb{N}}$:
\[
g_0 \equiv 0\,,
\]
while $h_0$ is the solution of the problem
\[
\dot{h}^{\sharp}_0 (t,j,x)= {h}^{\sharp}_0 (t,-j,x)\,,\qquad j=1,-1\,,\qquad h_0 |_{t=0}=f_0\,,
\]
where $f_0 \in X_1^+$ is the initial datum.

Moreover,
\[
\dot{g}^{\sharp}_1 + {g}^{\sharp}_1 = 0 \,,\qquad g_1 |_{t=0}=f_0\,,
\]
\[
\dot{g}^{\sharp}_n + \mathbb{P}^{\sharp} (h_{n-1} ,g_{n-1} ) g^{\sharp}_n =
\mathbb{Q}^{\sharp} (g_{n-1} ,h_{n-1} ) \,,\qquad g_n |_{t=0}=f_0\,,
\]
for $n\geq 2$, and
\[
\dot{h}^{\sharp}_n  + \mathbb{P}^{\sharp} (g_{n-1} ,h_{n-1} ) h^{\sharp}_n =
\mathbb{Q}^{\sharp} (h_{n-1} ,g_{n-1} )  \,,\qquad h_n |_{t=0}=f_0\,,
\]
for $n\geq 1$.

Therefore we have
\[
g_1 (t,j,x) = f_0 (x-jt)e^{-t}\,
\]
\[
\| h_0 (t) \|_1 \leq \| f_0 \|_1 e^t\,,
\]
and
\[
0\leq g_1 (t,j,x) \leq f_0 (x-jt) \leq h_1 (t,j,x) \leq h_0 (t,j,x)\,.
\]

Assuming that
\begin{equation}\label{deltunia}
f_0 (x) \geq \mu \qquad\mathrm{a.a.}\; x\in\mathbb{T}
\end{equation}
where $\mu$ is a positive ($>0$) number, we obtain that on any compact interval $[0,T]$
\[
0 < \mu e^{-T} \leq g_1 \leq g_2 \leq g_3 \leq \ldots \leq h_2 \leq h_1 \leq h_0\,.
\]

Therefore, by the theorem of Beppo--Levy both sequences converge in $X_1$
\[
g_n \uparrow g\,,\qquad h_n \downarrow h\,,
\]
in $X_1$ and
\[
g\leq h\,.
\]

On any fixed time interval $[0,T]$ we may apply the inequality
\[
0\leq h^{\gamma} - g^{\gamma} \leq \gamma (\mu e^{-T})^{\gamma -1} (h-g)
\]
to show that $g=h$.

Thus $f:=g=h$ is a mild solution to the Eq. (\ref{2.3}) on $[0,T]$, $T>0$, in the following sense
\[
f\in L^{\infty} (0,T;X_1 )\,,
\]
\[
\dot{f}^{\sharp} (t)  + \int\limits_0^t P^{\sharp} (f,f) \, f^{\sharp}_n  =
f_0 + \int\limits_0^t Q^{\sharp} (f,f) \,.
\]
for a.a. $t\in [0,T]$. The similar argument as the one used for $g=h$ shows that the solution is
unique. Moreover, the solution is in $X_1^+$.

Therefore we have

\begin{theorem}\label{twociagach}
Let $f_0 \in X_1^+$ be such that (\ref{deltunia})is satisfied for  some positive $\mu$. Then, for any given
$T>0$ there exists a unique mild solution $f=f(t)\in X_1$ of Eq. (\ref{2.3}) with the initial datum $f_0$
on $[0,T]$. Moreover,

\begin{itemize}

\item $f(t)\in X_1^+$ for a.a. $t\in [0,T]$,

\item $\mu e^{-T} \leq f(t)$ for a.a. $t\in [0,T]$,

\item $\| f(t) \|_1 \leq \|f_0 \|_1 e^T$ for a.a. $t\in [0,T]$.

\end{itemize}

\end{theorem}

\begin{remark}{\em
Note that if (\ref{deltunia}) is satisfied, then $f$ is a mild solution of Eq.\ (\ref{2.3}) in which each term
is on the  right hand side is separated from zero and thus Fr\`{e}chet differentiable. Hence, the solution becomes
a classical solution provided the initial datum is differentiable. }
\end{remark}

\bigskip
\section{Formal macroscopic analysis}\label{formal}

\medskip
In this section we follow the strategy developed in the case of the so--called
{\sl hydrodynamic limits} for kinetic equations in \cite{La}.

It is well known that the macroscopic limit is defined by $\varepsilon\to 0$,
where $\varepsilon >0$ is a suitable small parameter (the {\sl Knudsen number}).
Classically, see e.g. \cite{La}, Equation (\ref{2.3}) can be considered in the
following two dimensionless forms
\begin{equation}\label{3.1}
\Big(\partial_t +j\partial_x\Big) f (t,j,x) = {\frac{1}{\varepsilon}}\,
Q[f] (t,j,x) \,,
\end{equation}
and
\begin{equation}\label{3.2}
\Big( \varepsilon\partial_t +j\partial_x\Big) f (t,j,x) = {\frac{1}{\varepsilon}}\,
Q[f] (t,j,x) \,,
\end{equation}
where the nonlinear operator $Q$ is given by the RHS of Eq. (\ref{2.3}). Equations (\ref{3.1}) and (\ref{3.2}) with the initial data
\begin{equation}\label{3.3}
f\Big|_{t=0} =F
\end{equation}
lead to singularly perturbed problems in the limit $\varepsilon\to 0$.

We introduce some basic notation. The $0$--th order
term in the Hilbert expansion is the function $M$ (a Maxwellian) which satisfies
\begin{equation}\label{3.4}
Q(M)=0\,.
\end{equation}
Therefore
\begin{equation}\label{3.5}
M(t,j,x)=\eta_j \,\varrho (t,x)\,,\qquad \eta_1 +\eta_{-1}=1\,,
\end{equation}
where either
\begin{equation}\label{3.6}
\eta_{1} =\eta_{-1} = {\frac{1}{2}}\,,
\end{equation}
or
\begin{equation}\label{3.7}
\eta_k =1\,,\qquad \eta_{-k} =0\,,
\end{equation}
for some $k\in\{ -1,1\}$, and $\varrho$ is the macroscopic density. Therefore
we conclude that two macroscopic limits can be considered for initial data close
to (\ref{3.5}) either with (\ref{3.6}) or (\ref{3.7}). This is, however, strictly related to stability
properties of the solutions.

Consider the first Fr{\'e}chet derivative
(linearization) $L$ of the operator $Q$ at $M$. It has the form
\begin{equation}\label{3.8}
Lf(j)={\frac{1-\gamma}{2}}\Big( f(-j)-f(j)\Big)
\end{equation}
if (\ref{3.6}) holds, and
\begin{equation}
Lf(j)=
\begin{cases}  \; f(-k) & \qquad{\mathrm{for}}\; j=k \;\\
-f(-k) & \qquad{\mathrm{for}}\; j=-k
\end{cases}
\label{23}
\end{equation}
if (\ref{3.7}) is satisfied. For simplicity the dependence of $t$ and $x$ is not explicitly
indicated.

The null-space of $L$, called the \textit{hydrodynamic space}, is spanned by the
Max\-wel\-lian $M$. Since the aim of asymptotic analysis is to isolate the slow dynamics in
the hydrodynamic space from fast, transient, behavior in the complementary space, the best
way is split the phase space into subspaces reducing $L$. Hence, we call the \textsl{kinetic space}
the space spanned by all remaining eigenvectors of $L$ (corresponding to eigenvalues with negative
real parts). The kinetic space will be denoted by $W$. This can be done by spectral projections.
For a general matrix $A$ having 0 as the dominant simple eigenvalue with the corresponding
eigenvector with $ e$ and $ e^*$ being the eigenvector of the transpose $A^*$ corresponding to 0
and normalized so as $e\cdot e^* =1$, the spectral projection $\mc P$  onto the null--space of
$A$ is given by
\[
\mc P y = (e^*\cdot y) e\,.
\]
Then clearly the complementary projection $\mc P^\perp $ is defined as $\mc P^\perp = I -\mc P$,
where $I$ is the identity matrix. Specifying this theory to the case in hand, we see that in the
diffusive case (\ref{3.8}) the null spaces of $L$ and $L^*=L$ are spanned by, respectively,
$(1/2,1/2)$ and $(1,1)$ (where we have taken into account the normalization).  Thus, in this case we have
\begin{equation}\label{3.10}
\mathcal{P} \left( f(1)\,,\, f(-1)\right)=\left( {\frac{f(1)+f(-1)}{2}}\,,\,
{\frac{f(1)+f(-1)}{2}}\right)
\end{equation}
and
\begin{equation}\label{3.11}
\mathcal{P}^{\perp} \left( f(1)\,,\, f(-1)\right)=\left( {\frac{f(1)-f(-1)}{2}}\,,\,
{\frac{f(-1)-f(1)}{2}}\right).
\end{equation}
The situation is, however, different in both cases described by (\ref{3.7}). For $k=1$ the eigenvector
of $L$ corresponding to 0 is given by $(1,0),$ whereas $(0,1)$ is the respective eigenvector for $k=-1$.
For the transpose matrix the eigenvector is given by $(1,1)$ in either case. Thus, we obtain
\begin{equation}\label{3.10a}
\mathcal{P} \Big( f(1)\,,\, f(-1)\Big)=\Big( {{f(1)+f(-1)}}\,,\,
0\Big)
\end{equation}
and
\begin{equation}\label{3.11a}
\mathcal{P}^{\perp} \Big( f(1)\,,\, f(-1)\Big)=\Big( -f(-1)\,,\, f(-1)\Big),
\end{equation}
if $k=1,$ and
\begin{equation}\label{3.10b}
\mathcal{P} \Big( f(1)\,,\, f(-1)\Big)=\Big( 0,\, {{f(1)+f(-1)}}
\Big)
\end{equation}
with
\begin{equation}\label{3.11b}
\mathcal{P}^{\perp} \Big( f(1)\,,\, f(-1)\Big)=\Big( f(1)\,,\, -f(1)\Big)
\end{equation}
for $k=-1$.
We observe that the above expression give a preliminary  justification for the terms diffusive and
aligned picture. Indeed, in the first case the total density $\varrho = f(1)+f(-1)$ is uniformly
spread between the two available directions of motion whereas in the other two the total density
is concentrated along one of the directions of motion.
Let us introduce an operator $D$ by the formula
\[
D (f(1,x),f(-1,x)) = (\p_t f(1,x)+\p_x f(1,x), \p_t f(-1,x)-\p_x f(-1,x))\,.
\]
Referring to the problem (\ref{3.1}) we note (see \cite{La}) that
the abstract zeroth order (`Euler') macroscopic approximation is given by
\begin{equation}\label{E}
\mathcal{P}DM=0\,,
\end{equation}
whereas the first (`Navier-Stokes') order macroscopic approximation is
given by
\begin{equation}\label{NS}
\mathcal{P}DM=-\varepsilon\mathcal{P}D \mc P^\perp L^{-1}{\mathcal{P}}^{\perp} DM
\end{equation}
in terms of the macroscopic density $\varrho$.

Specifying for our three cases, we obtain the following formulae
\[
DM = \left\{\begin{array}{lcl} \frac{1}{2}(\p_t\varrho
+ \p_x\varrho, \p_t\varrho - \p_x\varrho)&{\rm in}& {\rm the\;diffusive\; case} \\
(\p_t\varrho + \p_x\varrho, 0)&{\rm in}& {\rm the\;aligned\; case,\;k=1}\\
 (0, \p_t\varrho - \p_x\varrho)&{\rm in}& {\rm the\;aligned\; case,\;k=-1},
\end{array}
\right.
\]
which, upon application of an appropriate $\mc P$,  yield the
following Euler approximating equations for the diffusive and aligned regimes with
$k=1,-1$, respectively,
\begin{equation}\label{hydapp1}
\begin{array}{ll}
\p_t\varrho & = 0\,, \\
\p_t\varrho +\p_x\varrho &  =  0\,, \\
\p_t\varrho -\p_x\varrho &  =  0\,.
\end{array}
\end{equation}
To derive the Navier-Stokes picture, first we find
\[
\mc P^\perp DM = \left\{\begin{array}{lcl} (\p_x\varrho, - \p_x\varrho)&
{\rm in}& {\rm the\;diffusive\; case,}\\
(0, 0)&{\rm in}& {\rm the\;aligned\; case,\;k=1,}\\
(0,0)&{\rm in}& {\rm the\;aligned\; case,\;k=-1},
\end{array}
\right.
\]
Clearly, the operator $L$ is not invertible but becomes such if reduced to the kinetic
space $W$ which, in the diffusive case, is spanned by $(1,-1)$. We see that $\mc P^\perp DM \in W$
and the solution to $L(y,z)= (\p_x \varrho, -\p_x\varrho)$ subject to $y=-z$ is
$(-\p_x\varrho, \p_x \varrho)/(1-\gamma)$. Next,
\[
D {\mathcal{P}}^{\perp}L^{-1}{\mathcal{P}}^{\perp} DM
= (1-\gamma)^{-1} (-\p_t\p_x\varrho -\p^2_x\varrho, \p_t\p_x\varrho -\p^2_x\varrho)
\]
and finally
\[
\mathcal{P}D \mc P^\perp L^{-1}{\mathcal{P}}^{\perp} DM= -(1-\gamma)^{-1} (\p^2_x\varrho, \p^2_x\varrho)
\]
which, combined with $\mc P D M = (\p_t\varrho, \p_t\varrho),$ yields the first order macroscopic
approximation (\ref{NS}) in the form
\begin{equation}\label{macro1}
\partial_t \varrho ={\frac{\varepsilon}{1-\gamma}}\partial_x^2 \varrho \,.
\end{equation}
Clearly, the Navier--Stokes approximation in either aligned case is again the respective Euler approximation.

Next we consider the diffusive picture regime  for problem the parabolic scaling  (\ref{3.2}).
Performing similar considerations we see that the zeroth order macroscopic approximation  is given by
\begin{equation}\label{macro2}
\partial_t \varrho ={\frac{1}{1-\gamma}}\partial_x^2 \varrho \,.
\end{equation}

Equations (\ref{macro1}) and (\ref{macro2}) are linear diffusion equations for $0<\gamma <1,$
with the diffusion coefficients ${\frac{\varepsilon}{1-\gamma}}$ and ${\frac{1}{1-\gamma}}$, respectively.
Linear dependence on $\varepsilon$ of the diffusion coefficient in the limit equation (\ref{macro1})
in the case of scaling (\ref{3.1}) is typical (cf. \cite{La, MBbook}). In the case of $\gamma >1$
Eqs. (\ref{macro1}) and (\ref{macro2}) are the (linear) backward diffusion equations. The backward diffusion
equation is less mathematically tractable than the diffusion equation --- see however \cite[Example 8.6]{dL}.

Thus, we cannot expect any stability properties in the diffusive regime  if $\gamma >1$. On the
other hand, one may expect that for $0< \gamma <1$ the asymptotic relationships
between the solutions of (\ref{3.1}) or (\ref{3.2}) and (\ref{macro1}) or
(\ref{macro2}), respectively, hold.

Next, let us consider the aligned regime for Eq. (\ref{3.1}) with
initial data (\ref{3.3}). Let $k\in\{-1,1\}$ be fixed and, to avoid technical difficulties we consider
only $\gamma >1$. We have seen that the macroscopic limit at both Euler and Navier--Stokes levels  does not depend on
the choice of $\gamma$ and is given by
\begin{equation}\label{3.14a}
\partial_t \varrho + k\,\partial_x \varrho = 0.
\end{equation}
Equation (\ref{3.14a}) is a linear first order equation and for any smooth initial datum  $F$ its
unique solution is given by
\begin{equation}\label{3.15}
\varrho (t,x)=F(k,x- k\, t)\,,\qquad t\geq 0\,,\quad x\in \mbb T\,.
\end{equation}
Clearly, the solution inherits all properties of $F$ such as boundedness, differentiability, etc.
This solution  is a traveling wave (in the direction $k$).

Next we investigate in more detail to aligned picture and, in particular, we provide a more refined
expansion formulae. To simplify notation, let us fix $k=1$; the calculation for $k=-1$ being symmetric.

Then, following (\ref{3.10a}) and (\ref{3.11a}), we have the hydrodynamic part of the solution
given by $\varrho = f_1+f_{-1}$ and the kinetic part defined as  $w = f_{-1}$ so that
\begin{equation}
(f_1,f_{-1}) =  (1,0)\varrho + ( -1,1)w.
\label{neww}
\end{equation}
Using $f_1=\varrho -w$ and $f_{-1}= w$, we see that the system (\ref{3.1}) can be written in hydrodynamic
and kinetic variables as
\begin{eqnarray}
\p_t \varrho + \p_x \varrho -2 \p_x w  &=& 0\,, \nn\\
\p_t w-\p_x w & =& \frac{1}{\e} \,\frac{(\varrho-w)w^\gamma -w(\varrho -w)^\gamma}{(\varrho-w)^\gamma +w^\gamma}\,.
\label{projeq}
\end{eqnarray}
First we look for the bulk approximation.

We perform the asymptotic expansion in the spirit of the Chapman--Enskog, that is, we do not expand
the hydrodynamic part and seek approximation $(\varrho, w) =(\varrho, w_0+\e w_1)$, see e.g. \cite{MBbook}. This produces
\begin{eqnarray}
&&\p_t \varrho + \p_x \varrho -2 \p_x (w_0+\e w_1) = 0\,, \nn\\
&&\e\p_t (w_0+\e w_1)-\e\p_x (w_0+\e w_1)\label{projeq1}\\
&&\phantom{xxx}= \frac{ \big( \varrho- ( w_0+\e w_1 ) \big) \big( w_0+\e w_1 \big)^\gamma - \big( w_0+\e w_1 \big) \big( \varrho
-(w_0+\e w_1) \big)^\gamma}{ \big( \varrho- ( w_0+\e w_1 ) \big)^\gamma + \big( w_0+\e w_1 \big)^\gamma}\,.\nn
\end{eqnarray}
Taking $O(1)$ terms on both sides, we obtain
\[
0 = \frac{(\varrho-w_0)w_0^\gamma -w_0(\varrho -w_0)^\gamma}{(\varrho-w_0)^\gamma +w_0^\gamma}
\]
which yields $w_0=0$ or $w_0 =\varrho$. If we recognize that the kinetic part for $k=1$ is the hydrodynamic
part for $k=-1$, then the latter option simply give the $k=-1$ case and we are left with $w_0=0$. To get $w_1$ we write
\begin{eqnarray*}
\frac{(\varrho-\e w_1)(\e w_1)^\gamma -\e w_1(\varrho -\e w_1)^\gamma}{(\varrho-\e w_1)^\gamma +(\e w_1)^\gamma} &=&
\frac{-\e\varrho^\gamma w_1 +O(\e^{\min\{2,\gamma\}})}{\varrho^\gamma(1-\e\gamma w_1\varrho^{-1} +O(\e^{\min\{2,\gamma\}}))}\\
&=& -\e w_1 + O(\e^{\min\{2,\gamma\}})
\end{eqnarray*}
and, since the remaining LHS terms in (\ref{projeq1}) are $O(\e^2),$ it follows that $w_1=0$. Hence, we recovered the
relevant, second, equation from (\ref{hydapp1}). We note that in this model and the presented level of expansion
Chapman--Enskog and Hilbert methods yield the same approximation.

\begin{remark}{\em
By induction we can prove that all terms of the expansion $w=w_0+\e w_1+\ldots +\e^n w_n +\ldots$ are zero.
Indeed, if $w_k=0$ for $k=0,1,\ldots,n-1$, the calculation as above gives
\begin{eqnarray*}
&&\frac{(\varrho-\e^n w_n)(\e^n w_n)^\gamma -\e^n w_n(\varrho -\e^n w_n)^\gamma}{(\varrho-\e^n w_n)^\gamma +(\e^n w_n)^\gamma}\\
&&\phantom{xx} =
\frac{-\e^n\varrho^\gamma w_n +O(\e^{\min\{2n,\gamma n\}})}{\varrho^\gamma(1-\e^n\gamma w_n\varrho^{-1} +O(\e^{\min\{2n,\gamma n\}}))}
= -\e^n w_n + O(\e^{\max \{ 2n,\gamma n \}})
\end{eqnarray*}
and  since the LHS of (\ref{projeq1}) now is $O(\e^{n+1}),$ we find $w_n=0$.}
\end{remark}

Next we incorporate the initial layer correction and look at the approximation
$(\varrho(t), w(t)) \approx (\bar\varrho(t) + r( \frac{t}{\varepsilon} ),  h({\frac{t}{\varepsilon}}))$,
where $\bar\varrho$ is the solution to the initial value problem
\begin{equation}
\p_t \bar\varrho + \p_x \bar\varrho = 0, \qquad \bar\varrho(0,x) = \varrho_0(x) = F(1,x).
\label{barrho}
\end{equation}
Then, denoting $\tau= \frac{t}{\varepsilon}$, the initial layer corrections $r$, $h$ are obtained by an asymptotic expansion of
\begin{eqnarray}
\p_t \bar\varrho +\e^{-1}\p_\tau r + \p_x (\bar\varrho +r) -2 \p_x h &=& 0 \,, \nn\\
\p_\tau h -\e\p_x h &=&  \frac{(\varrho+r-h)  h^\gamma -h(\varrho+r -h)^\gamma}{(\varrho+r-h)^\gamma +h^\gamma}\,.
\label{projeq2}
\end{eqnarray}
Since we expect $r$ to decay exponentially to 0 as $\tau\to \infty$, from the first equation we obtain $r=0$.
Then, setting $\e=0$ in the second equation and taking into account
that $\bar\varrho(t,x)|_{\e=0} = \bar\varrho(\e \tau)|_{\e=0} =  \varrho_0$ we obtain
\begin{eqnarray}
\partial_{\tau}h & = &\frac{ (\varrho_0 -h) h^{\gamma}- h (\varrho_0 -h)^{\gamma} }{(\varrho_0 -h)^{\gamma}+h^{\gamma}},\label{ileq}\\
 h (0,x) & = & F(-1,x)\,.\label{ileq2}
\end{eqnarray}

\bigskip
\section{Aligned picture}\label{aligned}

We consider the `{\textsl{aligned picture}}' for Eq. (\ref{3.1}) in the case $\gamma >1$.
We assume that $k\in\{ -1,1 \}$ is fixed.

Let $x\in{\mathbb{T}}$, where ${\mathbb{T}}$ is the $1$--dimensional torus. All functions are then interpreted
as  periodic functions on the unit interval $[0,1]$.

Let $X_{\infty}$ be the Banach space of continuous functions defined on ${\mathbb{T}}$ with the norm
\[
\| f\|_{\infty} = \sum\limits_{j\in\{ -1,1\}}\;\sup\limits_{x\in{\mathbb{T}}}\; |f(j,x)|\,.
\]

We consider Eq. (\ref{3.1}), i.e.
\begin{equation}\label{4.1}
\varepsilon \Big( \partial_t f_j  + j \,\partial_x f_j \Big) =
\frac{ \chi ( f_j + f_{-j} >0)}{f_j^{\gamma} + f_{-j}^{\gamma}}
\Big(  f_{-j} f_j^{\gamma} - f_j f_{-j}^{\gamma} \Big)\,, \qquad
j=k,-k\,,
\end{equation}
where $f_j (t,x)=f(t,j,x)$.

\medskip
Following the formal results derived in the previous section, we are looking for the solution
to Eq. (\ref{4.1}) in the following form
\begin{equation}\label{muzumba}
f_k (t,x) = \bar\varrho(t,x) - h \left( \frac{t}{\varepsilon} ,x\right) + \e u(t,x)\,, \quad f_{-k} (t,x)
= h \left( \frac{t}{\varepsilon} ,x\right) + \e v(t,x)\,,
\end{equation}
where $(\e u,\e v)$ is the error of the approximation. We recall that $\varrho =\varrho(t,x)$
is the `{\textsl{bulk solution}}' , $h=h (\tau ,x)$
is the `{\textsl{initial layer solution}}' term, $\tau=\frac{t}{\varepsilon}$ is the stretched time variable.
Then the error satisfies
\begin{eqnarray}
&&\varepsilon^2 \big( \partial_t u  + k \,\partial_x u \big) - \partial_{\tau} h - \varepsilon k\partial_x h \nn\\
&&\phantom{x}=\frac{ \chi \big(\bar \varrho + \varepsilon (u +v)\! >\!0 \big) \Big( \big( h+ \varepsilon v \big)
\big( \bar\varrho - h + \varepsilon u \big)^{\gamma}\!
- \! \big( \bar\varrho -h + \varepsilon u \big) \big( h+ \varepsilon v \big)^{\gamma} \Big)}{ \big( \bar\varrho
- h + \varepsilon u \big)^{\gamma} + \big( h+ \varepsilon v \big)^{\gamma}} \,, \nn\\
&&\varepsilon^2 \big( \partial_t v  - k \,\partial_x v \big) + \partial_{\tau} h - \varepsilon k\partial_x h \nn\\&&\phantom{x} =
\frac{ \chi \big( \varrho + \varepsilon (u + v)\! >\!0 \big) \Big(  \big( \bar\varrho -h
+ \varepsilon u \big) \big( h+ \varepsilon v \big)^{\gamma}\!
- \! \big( h+ \varepsilon v \big) \big( \bar\varrho -h
+ \varepsilon u \big)^{\gamma} \Big)}{ \big( \bar\varrho -h+ \varepsilon u \big)^{\gamma}
+ \big( h+ \varepsilon v \big)^{\gamma}}\,,\nn\\
\label{4.2}
\end{eqnarray}
where $\bar\varrho$ satisfies (\ref{3.15}), $h=h(\tau ,x)$ satisfies (\ref{ileq}), (\ref{ileq2})
and $\varrho_0 (x)= \varrho (0,x)=F(k,x)$.

First we consider the initial layer problem for $h$. We have

\begin{lemma}\label{lemat}
Let $k\in\{ -1,1\}$ and the initial data $F\in X$ be such that
\begin{equation}\label{warunek3}
F(-k,x)\geq\mu \,, \qquad\forall\; x\in\mathbb{T}\,,
\end{equation}
for some $\mu >0$. There exists $c_\gamma>1$ such that, \textbf{if}
\begin{equation}\label{warunek}
F(k,x)> c_\gamma\, F(-k,x) \qquad \forall\; x\in{\mathbb{T}}\,,
\end{equation}
\textbf{then} there exists a unique solution $h=h(\tau ,x)$ to the problem (\ref{ileq}), (\ref{ileq2})
for any $\tau >0$ and
\begin{equation}\label{zanik}
0< h(\tau ,x) \leq  F (-k,x) e^{-\delta\tau} \qquad\forall\; \tau >0 \quad \forall\;
x\in\mathbb{T}\,,
\end{equation}
for some $\delta >0$. Moreover  if, additionally, $F$ is continuously differentiable and $c_\gamma$ is sufficiently large,
then the solution $h=h(\tau ,x)$ to  (\ref{ileq}), (\ref{ileq2}) is continuously
differentiable and
\begin{equation}\label{zanik2}
|\partial_x h(\tau ,x)| \leq \mathrm{const.} \Big( \| F \|_{\infty} + \| \partial_x F \|_{\infty}
\Big) e^{-\delta_1 \tau} \qquad\forall\; \tau >0 \quad \forall\; x\in\mathbb{T}\,,
\end{equation}
for some $\delta_1 >0$.
\end{lemma}

\medskip
\begin{proof}
Problem (\ref{ileq}), (\ref{ileq2}) is an ODE problem, where the variable $x$ is a parameter. For $\gamma >1$ the
operator defined by the RHS of (\ref{ileq}) is Lipschitz continuous on
$C(\mathbb{T})$. Under the condition (\ref{warunek}) for a.a. $x\in\mathbb{T}$
and $\tau =0$ the RHS of Eq. (\ref{ileq}) is negative. Therefore $h$ strictly decreases with
respect to the variable $\tau$ for all $x\in\mathbb{T}$. Now,
\begin{eqnarray*}
&&\frac{1}{(\varrho_0 -h)^{\gamma}+h^{\gamma}}
\Big( (\varrho_0 -h) h^{\gamma}- h (\varrho_0 -h)^{\gamma} \Big) \\
&&\phantom{xxxxxxxx}= -h \frac{(\varrho_0-h)((\varrho_0-h)^{\gamma -1}-h^{\gamma-1})}{(\varrho_0 -h)^{\gamma}+h^{\gamma}}.
\end{eqnarray*}
Since, by (\ref{warunek}),
$$
0\leq \vartheta:=\frac{h}{\varrho_0-h} \leq \frac{F(-k)}{F(k)-F(-k)}\leq \frac{1}{c_\gamma-1},
$$
we can write
\begin{equation}
\frac{(\varrho_0-h)((\varrho_0-h)^{\gamma -1}- h^{\gamma-1})}{(\varrho_0 -h)^{\gamma}+h^{\gamma}} =
\frac{1-\vartheta^{\gamma-1}}{1+\vartheta^{\gamma}}:=\delta>0
\label{55}
\end{equation}
uniformly in $x\in [0,1]$ and $\tau>0$. We note that $\delta$ can be made as close to 1 as we wish
by making $c_\gamma$ sufficiently large but since close to the equilibrium $h=0$ the left hand side
of (\ref{55}) is close to 1, $\delta$ could not exceed 1. This leads to the conclusion that
$h$ decays exponentially to $0$ for all $x\in\mathbb{T}$ and the rate of convergence can be controlled
by $e^{-\delta \tau}$ uniformly in $x\in\mathbb{T}$. Moreover it is strictly positive. Thus the unique solution of Eq.
(\ref{ileq}) satisfies (\ref{zanik}).

It is evident that if $F$ is continuously differentiable with respect to $x$, then also $h$ is.
We denote $h' =\partial_x h$ and $\varrho_0'=\partial_x \varrho_0$.
The function $h' (\tau, x)$ satisfies
\begin{eqnarray}
&&\partial_{\tau} h' =
\frac{\gamma (\varrho_0 -h)^{\gamma -1}(\varrho_0' -h')
+\gamma h^{\gamma -1}h'}{((\varrho_0 -h)^{\gamma}+h^{\gamma})^2} \,
(\varrho_0 -h)\, h\, \Big( (\varrho_0 -h)^{\gamma -1} -h^{\gamma -1}\Big) \nn\\
&&\phantom{xx}+ \frac{-h'(\varrho_0 -h)\big( (\varrho_0 -h)^{\gamma -1} -h^{\gamma -1} \big)
-h(\varrho_0' -h')\big( (\varrho_0 -h)^{\gamma -1}
-h^{\gamma -1} \big)}{(\varrho_0 -h)^{\gamma}+h^{\gamma}} \nn\\
&&\phantom{xx}- (\gamma -1) \frac{h (\varrho_0 -h)^{\gamma -1} (\varrho_0' -h')
-h^{\gamma -1} (\varrho_0 -h) h'}{(\varrho_0 -h)^{\gamma}+h^{\gamma}}\,.
\label{mubmub}
\end{eqnarray}
The first term on the RHS of (\ref{mubmub}) does not cause any difficulty: the condition (\ref{warunek})
is sufficient to ensure the negativity of the term multiplying $h'$. Here we are using (\ref{zanik}).
Considering the second and the third terms on the RHS we see that if
\begin{eqnarray}
&&-(\varrho_0 -2h)\Big( (\varrho_0 -h)^{\gamma -1} -h^{\gamma -1} \Big) + (\gamma -1)h
(\varrho_0 -h)^{\gamma -1} \nn\\
&&\phantom{xxxxxxxxxxxxxxxxxxxxxxxxxx}+  (\gamma -1)h^{\gamma -1} (\varrho_0 -h) <0
\label{bumbum}
\end{eqnarray}
for all $x\in\mathbb{T},$ then (\ref{zanik2}) is satisfied. It is easy to see that
if  $c_{\gamma}$ is sufficiently large (for example
$c_2\geq  {6}/{(3-\sqrt{3})}$), then
(\ref{bumbum}) follows. This ends the proof.
\end{proof}

\bigskip

We may state the following result

\begin{theorem}\label{th1}
Let the initial data $F\in X_{\infty}$ be  nonnegative functions with continuous second derivatives and
such that (\ref{warunek3}), for some $\mu >0$, is satisfied.
For any $T>0$ there exists $\varepsilon_0 >0$ and $c_{\gamma} >0$ such that \textbf{if}
\begin{equation}\label{warunek4}
\min\limits_{x\in\mathbb{T}} F(k,x) \geq c_{\gamma}\, \max\limits_{x\in\mathbb{T}} F(-k,x)
\end{equation}
is satisfied for given $k\in\{ -1,1\}$, \textbf{then} for
$\,\varepsilon\in \; [\, 0,\varepsilon_0 \, [\,$ the Cauchy Problem for Eq. (\ref{4.1}),
with the initial datum $F$, has a mild solution $f=f(t)$ in $X_{\infty}$ on $[0,T]$.
Moreover,
\begin{equation}
\sup\limits_{t\in [0,T]} \, \Big\| f(t) - \bor (t) - \mb{h} \left( \frac{t}{\varepsilon} \right)
\,\Big\|_{\infty} \leq c_T \varepsilon\,,
\end{equation}
where
\[
\begin{array}{lll}
\bor (t,k,x) & = & F(k,x-kt)\,, \\
\bor (t,-k,x) & = & 0\,, \\
\mb{h} \Big( \frac{t}{\varepsilon} ,k,x \Big) & = & -h\Big( \frac{t}{\varepsilon} ,x \Big)\,, \\
\mb{h} \Big( \frac{t}{\varepsilon} ,-k,x \Big) & = & h\Big( \frac{t}{\varepsilon} ,x \Big)\,,
\end{array}
\]
$\mb{h}=\bar{h}\Big( \frac{t}{\varepsilon} ,j,x \Big)$, $j=k,-k$,
\[
\mb{h} (0,x) =F(-k,x) \qquad\forall\; x\in{\mathbb{T}}\,.
\]
\end{theorem}

\medskip
\begin{proof}
Let $T>0$ be fixed. The proof follows by the analysis of the system (\ref{4.2}) that, by
(\ref{muzumba}), is
equivalent to Eq. (\ref{4.1}). The methods of Theorem \ref{twociagach} may be used to
show the existence and uniqueness of the mild solution $f=f(t)$ in $X_{\infty}$ for any fixed
$\varepsilon >0$. The solution satisfies
\begin{itemize}

\item $f(t)\in X_{\infty}$ and $f(t)\geq 0$ for a.a. $t\in [0,T]$,

\item $\mu e^{-\frac{T}{\varepsilon}} \leq f(t)$ for a.a. $t\in [0,T]$,

\item $\| f(t) \|_{\infty} \leq \|f_0 \|_{\infty} e^{\frac{T}{\varepsilon}}$ for a.a. $t\in [0,T]$.
\end{itemize}

We may rewrite Eq. (\ref{4.1}) in the following form
\begin{equation}\label{4.3}
\varepsilon \partial_t f^{\sharp}_{j}  = \Bigg(
\frac{ \chi ( f_j + f_{-j} >0) \, f_{-j} f_j }{f_j^{\gamma} + f_{-j}^{\gamma}}
\Big(  f_j^{\gamma -1} - f_{-j}^{\gamma -1} \Big) \Bigg)^{\sharp}\,, \qquad
j=k,-k\,,
\end{equation}
where "$\sharp$" is given by (\ref{szarp}). Under the assumption (\ref{warunek3}) the
functions $f(t,k,x+kt)$ and $f(t,-k,x-kt)$ are increasing and decreasing functions
of $t$ (for any fixed $x$), respectively, and we have
\begin{equation}\label{4.4}
\begin{array}{lll}
& & f(t,k,x) \;\geq\; F(k,x-kt) \;\geq\; \min\limits_{y\in\mathbb{T}} F(k,y) \; > \\
& > & c_{\gamma} \max\limits_{y\in\mathbb{T}} F(-k,y) \;\geq\; c_{\gamma} F(-k,x+kt)
\;\geq\; c_{\gamma} f(t,-k,x)\,,
\end{array}
\end{equation}
for any $t>0$ and any $x\in\mathbb{T}$.

Now we assume (\ref{muzumba}) and consider Eq. (\ref{4.2}), taking in consideration that
$\bar\varrho + \varepsilon u +\varepsilon v >0$, we obtain
\begin{eqnarray}
\partial_t u  + k \,\partial_x u & \!\!=\!\! & \frac{1}{\varepsilon^2} \, \partial_{\tau} h
+ \frac{1}{\varepsilon} \, k\partial_x h \nn\\
& &\phantom{xx}+
\frac{\big( h+ \varepsilon v \big)\big( \bar\varrho - h + \varepsilon u \big)^{\gamma}
- \big( \bar\varrho -h + \varepsilon u \big) \big( h+ \varepsilon v \big)^{\gamma}}{ {\varepsilon^2}
\, \Big( \big( \bar\varrho - h + \varepsilon u \big)^{\gamma}
+ \big( h+ \varepsilon v \big)^{\gamma} \Big)}  \,, \nn\\
\partial_t v  - k \,\partial_x v & \!\!= \!\!& - \frac{1}{\varepsilon^2} \, \partial_{\tau} h
+ \frac{1}{\varepsilon} \, k\partial_x h \nn\\
&& \phantom{xx}+
\frac{\big( \bar\varrho -h+ \varepsilon u \big) \big( h+ \varepsilon v \big)^{\gamma}
- \big( h+ \varepsilon v \big) \big( \bar\varrho -h+ \varepsilon u \big)^{\gamma}}{ {\varepsilon^2} \; \Big(
\big( \bar\varrho -h+ \varepsilon u \big)^{\gamma} + \big( h+ \varepsilon v \big)^{\gamma}\Big)}  \,.
\label{4.5}
\end{eqnarray}
Therefore, using (\ref{ileq}) and (\ref{ba-ba}),
\begin{eqnarray}
&&\partial_t u  + k \,\partial_x u =\nn\\
&&\phantom{xxx}\frac{1}{\varepsilon} k\partial_x h + \frac{1}{{\varepsilon^2}}
\Bigg( \frac{\big( \bar\varrho - h + \varepsilon u \big)^{\gamma}
\big( \bar\varrho +\varepsilon u +\varepsilon v \big)}
{\big( \bar\varrho - h + \varepsilon u \big)^{\gamma} + \big( h+ \varepsilon v \big)^{\gamma}}
- \frac{ \big( \bar\varrho - h \big)^{\gamma} \bar\varrho}{\big( \bar\varrho - h \big)^{\gamma}
+ h^{\gamma}} \Bigg) \nn\\
&&\phantom{xxx}-\frac{1}{\varepsilon} u
+ \frac{1}{{\varepsilon^2}} \Bigg( \frac{ h \big( \bar\varrho - h \big)^{\gamma}
- (\bar\varrho -h) h^{\gamma}}
{ \big( \bar\varrho - h \big)^{\gamma} + h^{\gamma}}
- \frac{ h \big( \varrho_0 - h \big)^{\gamma} - \big( \varrho_0 -h \big)
h^{\gamma}}{\big( \bar\varrho_0 - h \big)^{\gamma} + h^{\gamma}} \Bigg) \,,\nn\\
&&\partial_t v  - k \,\partial_x v = \nn\\
&&\phantom{xxx}\frac{1}{\varepsilon} k\partial_x h +  \frac{1}{{\varepsilon^2}}\Bigg(
\frac{ \big( h+ \varepsilon v \big)^{\gamma} \big( \bar\varrho +\varepsilon u + \varepsilon v \big)}{
\big( \bar\varrho -h+ \varepsilon u \big)^{\gamma} + \big( h+ \varepsilon v \big)^{\gamma}}
- \frac{ h^{\gamma} \bar\varrho}{\big( \bar\varrho -h \big)^{\gamma} + h^{\gamma} } \Bigg) \nn\\
&&\phantom{xxx}-\frac{1}{\varepsilon} v + \frac{1}{\varepsilon^2} \Bigg( \frac{ \big( \bar\varrho -h \big) h^{\gamma}
- h \big( \bar\varrho -h\big)^{\gamma}}{ \big( \bar\varrho -h \big)^{\gamma} + h^{\gamma}} - \frac{ \big( \varrho_0 -h\big) h^{\gamma}
- h \big( \varrho_0 -h \big)^{\gamma}}{ \big( \varrho_0 -h \big)^{\gamma} + h^{\gamma}} \Bigg) \,.\nn\\
\label{4.6}
\end{eqnarray}
Assume now that $\varrho$, $h$, $u$, $v$ are given and consider the following functions
\[
\Xi_1 (\varepsilon )=
\frac{\big( \bar\varrho - h + \varepsilon u \big)^{\gamma} \big(\bar \varrho + \varepsilon u  + \varepsilon v \big)}
{\big( \bar\varrho - h + \varepsilon u \big)^{\gamma} + \big( h+ \varepsilon v \big)^{\gamma}}\,,
\]
and
\[
\Xi_2 (\varepsilon )=
\frac{\big( h + \varepsilon v \big)^{\gamma} \big( \bar\varrho + \varepsilon u  + \varepsilon v \big)}
{\big( \bar\varrho - h + \varepsilon u \big)^{\gamma} + \big( h+ \varepsilon v \big)^{\gamma}}\,,
\]
We have
\[
\Xi_i (\varepsilon ) - \Xi_i (0) = \varepsilon\, \int\limits_0^1 \,\Xi'_i (\theta\varepsilon ) \mathrm{d} \theta
\,,\quad i=1,2\,.
\]

The derivatives $\Xi'_i (\varepsilon )$, $i=1,2$, may be written as the sum
\[
\Xi'_i (\varepsilon ) = \Gamma_{i,1}  (\varepsilon ) u + \Gamma_{i,2} (\varepsilon ) v \,,
\]
where, denoting $\Theta =  \bar\varrho -h +\varepsilon u $ and $\Lambda = h+\e v,$
\begin{eqnarray*}
\Gamma_{1,1} (\varepsilon )& = & \frac{ \Theta^{2\gamma}
+ \gamma \Theta^{\gamma -1} \Lambda^{\gamma}
\big( \bar\varrho + \varepsilon u + \varepsilon v \big) + \Theta^{\gamma}
\Lambda^{\gamma} }{\Big( \Theta^{\gamma}
+ \Lambda^{\gamma} \Big)^2} \,, \\
\Gamma_{1,2} (\varepsilon )& = & \frac{\Theta^{2\gamma}
- \gamma \Theta^{\gamma} \Lambda^{\gamma -1}
\big( \bar\varrho + \varepsilon u + \varepsilon v \big) + \Theta^{\gamma}
\Lambda^{\gamma}}{\Big(\Theta^{\gamma}
+ \Lambda^{\gamma} \Big)^2}\,, \\
\Gamma_{2,1} (\varepsilon )& = & \frac{\Theta^{\gamma}
\Lambda^{\gamma} - \gamma \Theta^{\gamma -1}
\Lambda^{\gamma} \big( \bar\varrho + \varepsilon u + \varepsilon v \big)
+ \Lambda^{2\gamma}}{\Big( \Theta^{\gamma}
+ \Lambda^{\gamma} \Big)^2} \,, \\
\Gamma_{2,2} (\varepsilon )& = & \frac{ \Theta^{\gamma}
\Lambda^{\gamma}
+ \gamma \Theta^{\gamma} \Lambda^{\gamma -1}
\big( \bar\varrho + \varepsilon u + \varepsilon v \big)
+ \Lambda^{2\gamma}}{\Big( \Theta^{\gamma}
+ \Lambda^{\gamma} \Big)^2} \,.
\end{eqnarray*}
Let us consider
\begin{eqnarray*}
\Gamma_{1,1} (\varepsilon ) &=& \frac{ 1+ \frac{\Lambda^\gamma}{\Theta^\gamma}\frac{\varrho
+ \varepsilon u + \varepsilon v}{\Theta}+ \frac{\Lambda^\gamma}{\Theta^\gamma} }
{\Big( 1+ \frac{\Lambda^\gamma}{\Theta^\gamma}\Big)^2}  \\
&=& \left(1 +(h+\e v) \frac{\Lambda^{\gamma-1}}{\Theta^\gamma}\frac{\varrho
+ \varepsilon u + \varepsilon v}{\Theta}
+ \frac{\Lambda^{\gamma-1}}{\Theta^\gamma}\right)\left( 1 - \frac{\Lambda^\gamma}{\Theta^\gamma}
+ \ldots\right)^2\,.
\end{eqnarray*}
Now, by (\ref{warunek3}), (\ref{warunek4}) and (\ref{4.4}) we get
\begin{equation}\label{minest}
\begin{array}{l}
\bar\varrho - h + \e u  \geq  \min\{\bar\varrho-h, f_k\}  \\
\phantom{xxx}\geq \min\{F(k)-F(-k), F(k)\}\geq F(-k)\min\{c_\gamma-1,c_\gamma\} \geq  \mu (c_\gamma-1)\,,
\end{array}
\end{equation}
and
\begin{equation}
\frac{\Lambda}{\Theta}=\frac{h+\e v}{\bar\rho -h +\e u} \leq \frac{\max\{F(-k), h\}}{(c_\gamma -1)F(-k)} \leq
\frac{1}{c_\gamma -1} \to 0, \qquad c_\gamma \to \infty\,.
\end{equation}
Performing analogous estimates for other $\Gamma_{i,j}$, we see that
\begin{eqnarray}
&&\Gamma_{1,1} (\varepsilon )  +  \Gamma_{2,1} (\varepsilon ) = 1 \,, \nn\\
&&\Gamma_{1,2} (\varepsilon )  +  \Gamma_{2,2} (\varepsilon ) = 1 \,, \nn\\
&&|\Gamma_{1,1} (\varepsilon )| \leq 1 + \sigma_{1,1},  \qquad  |\Gamma_{2,1} (\varepsilon )| \leq \sigma_{2,1} \,, \nn\\
&&|\Gamma_{1,2} (\varepsilon )| \leq 1 - \sigma_{1,2},  \qquad  |\Gamma_{2,2} (\varepsilon )| \leq \sigma_{2,2} \,,
\label{4.7}
\end{eqnarray}
where the positive constants $\sigma_{i,j}$, $i,j=1,2$, may be as small as we want
if $c_{\gamma}$ is sufficiently large and also
\begin{eqnarray}
|\Gamma_{1,1} (\varepsilon )| & \leq & 1 + \sigma'_{1,1} \big( |h|+\varepsilon |v| \big)\,, \nn\\
|\Gamma_{2,1} (\varepsilon )| & \leq & \sigma'_{2,1} \big( |h|+\varepsilon |v| \big)\,, \nn\\
|\Gamma_{2,2} (\varepsilon )| & \leq & \sigma'_{2,2} \big( |h|+\varepsilon |v| \big)\,,
\label{4.8}
\end{eqnarray}
where again the positive constants $\sigma'_{1,1}$, $\sigma'_{2,1}$, $\sigma'_{2,2}$  may
be chosen  as small as we want by taking sufficiently large $c_\gamma$.

In order to estimate the last term of Eq. (\ref{4.6}) we assume that $\tau$ is fixed and consider
\[
\Psi (t) = \frac{(\varrho (t) -h (\tau )) h^{\gamma} (\tau ) - h (\tau ) (\varrho (t) -h (\tau ))^{\gamma}}
{( \varrho (t) -h (\tau ) )^{\gamma} + h^{\gamma} (\tau )\big)}\,.
\]
By (\ref{ba-ba}),
\[
\Psi(t) = \frac{h^\gamma\bar\varrho(t) }{(\bar\varrho(t)-h)^\gamma+h^\gamma}-h\,,
\]
so that
\[
\Psi(t) -\Psi(0) = \frac{h^\gamma\bar\varrho(t) }{(\bar\varrho(t)-h)^\gamma+h^\gamma}-
\frac{h^\gamma\varrho_0 }{(\varrho_0-h)^\gamma+h^\gamma}\,,
\]
and thus, using differentiability of the initial conditions, (\ref{4.4}) and (\ref{minest}),
\begin{eqnarray}
|\Psi(t) -\Psi(0)| &=& t h^\gamma |\rho_0'(\theta)| \left|\frac{(\bar\varrho(\theta)-h)^\gamma+h^\gamma -
\gamma\bar\varrho(\theta)(\bar\varrho(\theta) -h)^{\gamma -1}}{((\bar\varrho(\theta)-h)^\gamma+h^\gamma)^2}\right| \nn\\
&\leq&  t h^\gamma \| F'\|_\infty\frac{(2+\gamma)\|F\|^\gamma_\infty}{\mu^2(c_\gamma -1)^2} = C_\gamma th^\gamma\,.
\label{hga}
\end{eqnarray}
Then we have
\begin{eqnarray}
\partial_t u\!\!\!\!  &+& \!\!\!\!k \,\partial_x u + \frac{1}{\varepsilon} u  \nn\\
&=& \frac{1}{\varepsilon} k\partial_x h
+ \frac{u}{\varepsilon} \int\limits_0^1 \Gamma_{1,1} (\theta \varepsilon) \, \mathrm{d} \theta
+ \frac{v}{\varepsilon} \int\limits_0^1 \Gamma_{1,2} (\theta \varepsilon) \, \mathrm{d} \theta \,
+ \frac{t}{{\varepsilon}^2} \int\limits_0^1 \Psi'(\theta t)\,\mathrm{d}\theta, \nn\\
\partial_t v \!\!\!\! &-&\!\!\!\! k \,\partial_x v \!+\!\! \frac{1}{\varepsilon} v  \nn\\
&=& \frac{1}{\varepsilon} k\partial_x h
+ \frac{u}{\varepsilon} \int\limits_0^1 \Gamma_{2,1} (\theta \varepsilon)
\, \mathrm{d} \theta + \frac{v}{\varepsilon} \int\limits_0^1 \Gamma_{2,2} (\theta_2 \varepsilon) \,\mathrm{d} \theta
\!-\!\! \frac{t}{{\varepsilon}^2} \int\limits_0^1 \Psi'(\theta t)\,\mathrm{d}\theta \,.
\label{4.9}
\end{eqnarray}
Integrating (\ref{4.9}) along characteristics and using (\ref{4.7}) and (\ref{4.8}), we obtain
\begin{eqnarray}
\| u(t)\|_{\infty} &\leq& \frac{C}{\varepsilon} \int\limits_0^t
\left( 1+ \frac{s}{\varepsilon} \right) e^{-\frac{s}{\varepsilon} \bar\delta} \,\mathrm{d}s +
\frac{\sigma'_{1,1}}{\varepsilon} \int\limits_0^t   e^{-\frac{s}{\varepsilon} \delta}  \| u(s)\|_{\infty}\,\mathrm{d}s \nn\\
&&\phantom{xx}+ \sigma'_{1,1}\!\! \int\limits_0^t \| u(s)\|_{\infty}  \| v(s)\|_{\infty} \,\mathrm{d}s
+ \frac{1-\sigma_{1,2}}{\varepsilon}\!\! \int\limits_0^t\!\! \| v(s)\|_{\infty}\,\mathrm{d}s,\label{bum1}\\
\| v(t)\|_{\infty} &\leq& \frac{C}{\varepsilon} \int\limits_0^t e^{-\frac{t-s}{\varepsilon}}
\left( 1+ \frac{s}{\varepsilon} \right) e^{-\frac{s}{\varepsilon} \bar\delta} \,\mathrm{d}s \nn\\
&&\phantom{xx}+ \frac{\sigma_{2,1}}{\varepsilon}\!\! \int\limits_0^t e^{-\frac{t-s}{\varepsilon}} \| u(s)\|_{\infty}\,\mathrm{d}s
+ \frac{\sigma_{2,2}}{\varepsilon} \!\!\int\limits_0^t \!\!e^{-\frac{t-s}{\varepsilon}} \| v(s)\|_{\infty}\,\mathrm{d}s,\label{bum2}
\end{eqnarray}
where $C$ denotes a generic constant independent of $\e$ and $\bar\delta < \min\{\gamma\delta, \delta_1\}>0$,
see Lemma \ref{lemat} and (\ref{hga}). For convenience of calculations, we   fix $\bar\delta <1.$

We can re-write (\ref{bum2}) as
\begin{equation}\label{bum3}
\| v (t) \|_{\infty} \leq C + \sigma_{2,1} \| u\|_{\infty ,t}
+ \sigma_{2,2} \| v\|_{\infty ,t}\,,\qquad t\in [0,T]\,,
\end{equation}
where
$\| u\|_{\infty ,t} = \sup\limits_{s\in [0,t]}\| u(s)\|_{\infty}$.
Taking $c_\gamma$ large enough for $\sigma_{2,2}$ to be smaller then 1, we obtain
\begin{equation}\label{bum4}
\| v\|_{\infty ,t} \leq C(1+ \sigma_{2,1}) \| u\|_{\infty ,t} \,,\qquad t\in [0,T]\,.
\end{equation}

We may now refine  the estimate (\ref{bum2}), using (\ref{4.8}) and (\ref{zanik}), with $\delta$ replaced by $\bar\delta$,
\begin{eqnarray}
\| v(t)\|_{\infty}\!\! \!\!&\leq&\!\! \!\! \frac{C}{\varepsilon} \int\limits_0^t \!\!e^{-\frac{t-s}{\varepsilon}}
\left( 1+ \frac{s}{\varepsilon} \right) e^{-\frac{s}{\varepsilon} \bar\delta} \,\mathrm{d}s \nn\\
&&+ \frac{\sigma'_{2,1}}{\varepsilon}\!\! \int\limits_0^t \!\!e^{-\frac{t-s}{\varepsilon}}
e^{-\frac{s}{\varepsilon}\bar \delta} \| u(s)\|_{\infty}\,\mathrm{d}s
+ \frac{\sigma'_{2,2}}{\varepsilon}\!\! \int\limits_0^t \!\!e^{-\frac{t-s}{\varepsilon}} e^{-\frac{s}{\varepsilon} \bar\delta}
\| v(s)\|_{\infty}\,\mathrm{d}s \nn\\
&&+ \sigma'_{2,1}\!\! \int\limits_0^t \!\!e^{-\frac{t-s}{\varepsilon}} \| u(s)\|_{\infty} \| v(s)\|_{\infty} \,\mathrm{d}s +
\sigma'_{2,2} \!\!\int\limits_0^t \!\!e^{-\frac{t-s}{\varepsilon}} \| v(s)\|_{\infty}^2 \,\mathrm{d}s .
\label{bum5}
\end{eqnarray}
Therefore  we find
\begin{eqnarray}
\frac{1}{\varepsilon} \int\limits_0^t \| v(s)\|_{\infty} \,\mathrm{d}s\!\! &\leq & \!\!C +
\sigma'_{2,1} \| u\|_{\infty ,t} + \sigma'_{2,2} \| v\|_{\infty ,t} \nn \\
&&+\,\sigma'_{2,1} (T+\varepsilon ) \| u\|_{\infty ,t} \| v\|_{\infty ,t} +
\sigma'_{2,2} (T+\varepsilon ) \| v\|_{\infty ,t}^2.
\label{bum6}
\end{eqnarray}
By (\ref{bum1}), (\ref{bum4}) and (\ref{bum6}) we obtain
\begin{eqnarray*}
\| u(t)\|_{\infty} &\leq& C + \frac{\sigma'_{1,1}}{\bar\delta}  \| u\|_{\infty,t} +\sigma'_{1,1} T C
 \| u\|_{\infty,t}  (1+  \sigma_{2,1} \| u\|_{\infty ,t}) \\
 &&+ ({1-\sigma_{1,2}}) \left(C +
\sigma'_{2,1} \| u\|_{\infty ,t} + C\sigma'_{2,2}(1+  \sigma_{2,1} \| u\|_{\infty ,t}) \right.\nn\\
&&\left.+\sigma'_{2,1} (T+\varepsilon )C \| u\|_{\infty ,t}(1+  \sigma_{2,1} \| u\|_{\infty ,t}) \right.\\
&&+\left.
\sigma'_{2,2} (T+\varepsilon ) C^2(1+  \sigma_{2,1} \| u\|_{\infty ,t})^2 \right)
\end{eqnarray*}
so that
\begin{equation}\label{bum7}
\| u\|_{\infty ,t} \leq C + \sigma^{(1)} \| u\|_{\infty ,t} + \sigma^{(2)}
\| u\|_{\infty ,t}^2 \,.
\end{equation}
for any $t\in [0,T]$, where $\sigma^{(1)}$ and $\sigma^{(2)}$ are given by
\begin{eqnarray*}
\sigma^{(1)} &=& \frac{\sigma'_{1,1}}{\bar\delta} +\sigma'_{1,1}TC
+ (1-\sigma_{1,2})(\sigma'_{2,1} + C\sigma_{2,1} \sigma'_{2,2}) \\&&\phantom{xxxxxxxxxx}+ C(T+\e_0)(\sigma'_{2,1} + 2C
\sigma'_{2,2} \sigma_{2,1})\,,
\\
\sigma^{(2)} &=& \sigma'_{1,1} \sigma_{2,1}TC + (\sigma'_{2,1}\sigma_{2,1}+ C\sigma'_{2,2}\sigma^2_{2,1} \big)C(T+\e_0) .
\end{eqnarray*}
and tend to 0 with $\sigma_{i,j}\to 0$ and $\sigma'_{i,j}\to 0$ (that is, for $c_\gamma\to \infty$).
Thanks to this  we may choose $c_\gamma$ large enough for   $\sigma^{(1)} <1$ to obtain
\begin{equation}\label{bum8}
\| u\|_{\infty ,t} \leq C(1 + \sigma^{(2)} \| u\|_{\infty ,t}^2) \,,
\end{equation}
for any $t\in [0,T]$, where $C$ is a constant.

Let $y:= \| u\|_{\infty ,t}$, $\alpha := C$ and $\beta := C \sigma^{(2)}$.
Then, by (\ref{bum8}), we have
\begin{equation}\label{bum9}
0\leq - y + \alpha + \beta y^2\,.
\end{equation}

Considering the function $h(y)=\alpha - y + \beta y^2$, we may choose $\alpha >0$ and $\beta >0$
($\alpha\beta < \frac{1}{4}$) such that
\begin{equation}\label{bum10}
h(y_{\mathrm{min}})<0\,,\qquad\mathrm{for}\quad  y_{\mathrm{min}} = \frac{1}{2\beta }\,.
\end{equation}

Moreover $y\big|_{t=0}=0$ and $y$ cannot experience jumps. Therefore, in particular, we have
\begin{equation}\label{bum11}
y \leq \frac{1}{2\beta}\,,
\end{equation}
and $\|u\|_{\infty ,T}$ is bounded.
Then (\ref{bum4}) ends the proof.
\end{proof}

\bigskip\bigskip
\begin{center}
{\textbf{Acknowledgments}}
\end{center}
M.L. acknowledges a financial support from the Polish Ministry of Science and Higher
Education under the grant No. N N201 362536. J. B. acknowledges partial support from No. N N201 362536,
partial support from the National Scientific Centre of Poland under the grant  No. N N201 021133
and partial support from NRF Grant FA2007030300001.

\end{document}